\documentclass[12pt]{article}
\usepackage{comment}
\usepackage{amsmath}
\usepackage{amsthm}
\usepackage[font=small]{caption}
\usepackage{tikz}
\newtheorem*{theorem}{Theorem}
\theoremstyle{remark}
\newtheorem{remark}{Remark}
\newtheorem*{ack}{Acknowledgments}
\numberwithin{equation}{section}
\def\C{\mathbf{C}}

\def\R{\mathbf{R}}
\def\Ima{\mathrm{Im}\, }
\def\im{\operatorname{Im}}
\def\re{\operatorname{Re}}
\def\dist{\operatorname{dist}}
\def\meas{\operatorname{meas}}
\author{Walter Bergweiler and Alexandre Eremenko\thanks{Supported by NSF grant
DMS-1361836.}}
\title{On the Bank--Laine conjecture}
\date{}
\begin{document}
\maketitle

\begin{abstract} We resolve a question of Bank and Laine on the zeros
of solutions of $w^{\prime\prime}+Aw=0$ where $A$ is an entire function of finite order.

AMS 2010 Subj. Class: 34A20, 30D15. 

\end{abstract}

\section{Introduction and result}
The asymptotic distribution of zeros of solutions of linear differential equations 
with polynomial coefficients is described quite precisely by asymptotic integration methods;
cf.\ \cite{Hille} and \cite[Chapter~8]{Hille2}.
While certain differential equations with transcendental coefficients such as the Mathieu equation
were considered early on, the first general results concerning the 
frequency of the zeros of the solutions of
\begin{equation}\label{DE}
w^{\prime\prime}+Aw=0
\end{equation}
with a transcendental entire function $A$
appear to be due to  Bank and Laine~\cite{BL0,BL}.

For an entire function~$f$,
denote by $\rho(f)$ the order and by $\lambda(f)$ the exponent of convergence
of the zeros of $f$.
If $A$ is a polynomial of degree~$n$, then $\rho(w)=1+n/2$ for every
solution $w$ of~\eqref{DE},
while $\rho(w)=\infty$ for every solution $w$ if $A$ is transcendental.

Let $w_1$ and $w_2$ be linearly independent solutions of \eqref{DE}.
Bank and Laine proved that if $A$ is transcendental and
$\rho(A)<\frac12$, then 
$$\max\{\lambda(w_1),\lambda(w_2)\}=\infty.$$
It was shown independently by Rossi \cite{Ros} and Shen \cite{Shen}
that this actually holds for $\rho(A)\leq \frac12$.
Bank and Laine also showed that in the case of non-integer $\rho(A)$ we always have
\begin{equation}
\label{maxexpo}
\max\{\lambda(w_1),\lambda(w_2)\}\geq \rho(A),
\end{equation}
and they gave examples of functions $A$ of integer order for which there are
solutions $w_1$ and $w_2$ both without zeros. 

A problem left open by their work -- which later became known as the
Bank--Laine conjecture -- is whether 
$\max\{\lambda(w_1),\lambda(w_2)\}=\infty$ whenever
$\rho(A)$ is not an integer.
This question has attracted considerable interest; see \cite{L} for a survey, as well as, e.g.,~\cite{Gundersen}, \cite{HK}
and~\cite[Chapter~5]{Laine}.

We answer this question by showing that the estimate \eqref{maxexpo}
is best possible for a dense set of orders in the interval $(1,\infty)$.
\begin{theorem}
Let $p$ and $q$ be odd integers. Then there exists an entire function $A$ of order
$$\rho(A)=1+\frac{\log^2 (p/q)}{4\pi^2}$$
for which the equation \eqref{DE} has two linearly independent solutions 
$w_1$ and $w_2$ such that $\lambda(w_1)=\rho(A)$ while $w_2$ has no zeros.
\end{theorem}
By an extension of the method it should be possible to achieve any prescribed order
$\rho(A)>1$; see Remark~\ref{rem2} at the end. 

If
$w_1$ and $w_2$ are linearly independent solutions of \eqref{DE},
then the Wronskian $W(w_1,w_2)=w_1w_2^\prime-w_1^\prime w_2$
is a non-zero constant.
The solutions are called normalized if $W(w_1,w_2)=1$.

It is well-known that the ratio $F=w_2/w_1$ satisfies the Schwarz differential
equation (see, for example \cite{Hille2}):
$$
\label{schwarz}
S[F]:=\frac{F^{\prime\prime\prime}}{F^\prime}-\frac{3}{2}\left(\frac{F^{\prime\prime}}{F^\prime}\right)^2=2A.
$$
These meromorphic functions $F$ are completely characterized
by a topological property: they are locally univalent.
More precisely, consider the equivalence relation on meromorphic functions
$F_1\sim F_2$ if $F_1=L\circ F_2$, where $L$ is a fractional linear transformation.
Then the map $F\mapsto S[F]$ is a bijection between the equivalence classes of locally univalent
meromorphic functions and all entire functions.

Normalized solutions $w_1,w_2$ are recovered from $F$ by the formulas
$$w_1^2=\frac{1}{F'},\quad w_2^2=\frac{F^2}{F'}.$$
So zeros of $F$ are zeros of $w_2$ and poles of $F$ are zeros of $w_1$.

A meromorphic function $F$ is locally univalent if and only
if $E=F/F'$ is an entire function with the property that $E(z)=0$ implies
$E'(z)\in\{-1,1\}$. Such entire functions $E$ are called {\em Bank--Laine functions}.
If $w_1$ and $w_2$ is a normalized system of solutions of (\ref{DE})
and $F=w_2/w_1$, then 
$$E=\frac{F}{F'}=w_1w_2.$$ 
The converse is also true: every Bank--Laine function
is the product of two linearly independent solutions of (\ref{DE}). 

It turns out that the Schwarzian derivative has the following
factorization:
$$2S[F]=B[F/F'],$$
where
\begin{equation}\label{B}
B[E]:=-2\frac{E^{\prime\prime}}{E}+\left(\frac{E^\prime}{E}\right)^2-\frac{1}{E^2}.
\end{equation}
Thus every Bank--Laine function $E$ is a product of two linearly independent solutions
of (\ref{DE}) with $4A=B[E]$, a fact discovered by Bank and Laine \cite{BL0,BL}.

A considerable part of the previous research related to the Bank--Laine conjecture has concentrated on
the study of Bank--Laine functions.
There are a number of papers where Bank--Laine functions of finite order with various other properties are 
constructed~\cite{Alo2009,DL,Fle2009,L1,Lan2001,Lan2002,Mey2007}.
In all examples constructed so far, for which the order could be determined,
it was an integer. In our construction we have 
$\rho(E)=\rho(A)$; see Remark~\ref{rem1}. Thus our theorem also yields the
first examples of Bank--Laine functions of finite non-integral 
order.

In the proof of our theorem we use the fact that the functions $F$ have a topological characterization.
Starting with two elementary locally univalent functions, we paste them together by a quasiconformal
surgery. The resulting function is locally univalent, and the asymptotics of $\log|F/F'|$ can be
explicitly computed. A different kind of quasiconformal surgery was used in~\cite{DL,L}.

\begin{ack}
We are grateful to Jim Langley for a very detailed reading of the manuscript and many helpful suggestions.
The first author also thanks him for an illuminating discussion on the problem in 2012.
We are also thankful to David Drasin, Gary Gundersen and Ilpo Laine for their comments.
\end{ack}

\section{Proof of the theorem}

For every integer $m\geq 0$ we consider the polynomial
$$P_m(z)=\sum_{k=0}^{2m}(-1)^k \frac{z^k}{k!}.$$
Then the entire function
$$g_m(z)=P_m(e^z)\exp{e^z}$$
satisfies 
$$g_m^\prime(z)=\left(P_m'(e^z)+P_m(e^z)\right)e^z\exp {e^z}=\frac{1}{(2m)!}\exp \left(e^z+(2m+1)z\right)$$
and thus it
has the following properties:
\begin{itemize}
\item[a)] $g_m^\prime(z)\neq 0$ for all $z\in\C$,
\item[b)] $g_m$ is increasing on $\R$, and satisfies
$g_m(x)\to 1$ as $x\to-\infty$ as well as
$g_m(x)\to+\infty$ as $x\to+\infty$. 
\end{itemize}

From now on, we fix two distinct non-negative integers $m$ and~$n$,
and will sometimes omit them from notation. 
Notice that $g_m$ and $g_n$
are locally univalent entire functions.
We are going to restrict
$g_m$ to the upper half-plane $H^+$ and $g_n$ to the lower half-plane
$H^-$, and then paste them together, using a quasiconformal surgery, producing an entire function $F$.
Then our Bank--Laine function will be $E=F/F'$ and thus $A=B[E]/4$ as in \eqref{B}.

It follows from b) that there exists an increasing diffeomorphism
$\phi\colon\R\to\R$ such that $g_m(x)=(g_n\circ\phi)(x)$ for $x\in\R$.  Let
$$
k=\frac{2m+1}{2n+1}.
$$
We show that the asymptotic behavior
of the diffeomorphism $\phi$ is the following:
\begin{equation}\label{a1}
\phi(x)=x+O(e^{-x/2}),\quad \phi^\prime(x)\to 1,\quad x\to+\infty,
\end{equation}
and
\begin{equation}\label{a2}
\phi(x)=kx+c+O(e^{-\delta |x|}),\quad \phi^\prime(x)\to k, \quad x\to-\infty,
\end{equation}
with 
$$c=\frac{1}{2n+1}\log
\frac{(2n+1)!}{(2m+1)!} 
\quad \text{and} \quad 
\delta= \frac12 \min\{1,k\}.$$

In order to prove \eqref{a1}, we note that 
$$\log g_m(x)=e^x+O(x)=e^x\left(1+O(xe^{-x})\right), \quad x\to+\infty.$$
The equation $g_m(x)=g_n(\phi(x))$ easily implies
that $\frac23 x\leq \phi(x)\leq 2x$ for large~$x$.
Thus we also have 
\begin{align*}
\log g_n(\phi(x))&=e^{\phi(x)}\left(1+O\!\left(\phi(x)e^{-\phi(x)}\right)\right)\\
&=
e^{\phi(x)}\left(1+O\!\left(xe^{-2x/3}\right)\right), \quad x\to+\infty.
\end{align*}
Combining the last two equations we obtain 
$$e^{\phi(x)-x}=1+O(xe^{-2x/3}), \quad x\to+\infty,$$
from which the first statement in \eqref{a1} easily follows.
For the second statement in \eqref{a1} we use
\begin{equation}
\label{formula}
\phi'=\frac{g_m^\prime}{g_m}\frac{g_n\circ\phi}{g_n^\prime\circ\phi},
\end{equation}
so that
\begin{align*}
\phi'(x)&=\frac{(2n)!}{(2m)!}
e^{(2m+1)x-(2n+1)\phi(x)}\frac{P_n(e^{\phi(x)})}{P_m(e^x)}\\
&\sim  e^{(2m+1)x-(2n+1)\phi(x)+2n\phi(x)-2mx}\\
&=e^{x-\phi(x)}=1+o(1), \quad x\to+\infty.
\end{align*}
In order to prove \eqref{a2} we note that 
$$P_m(w)=e^{-w}+\frac{w^{2m+1}}{(2m+1)!}+O(w^{2m+2}),\quad w\to 0,$$
and thus 
$$P_m(w)e^w=1+\frac{w^{2m+1}}{(2m+1)!}+O(w^{2m+2}),\quad w\to 0.$$
Hence 
\begin{align*}
g_m(x)
&=1+\frac{e^{(2m+1)x}}{(2m+1)!}+O(e^{(2m+2)x})\\
&=1+\frac{e^{(2m+1)x}}{(2m+1)!}(1+O(e^x)),\quad x\to -\infty.
\end{align*}
The equation $g_m(x)=g_n(\phi(x))$ now yields
$$\frac{(2m+1)!}{(2n+1)!}e^{(2n+1)\phi(x)-(2m+1)x}=1+O(e^x)+O(e^{\phi(x)}),\quad x\to -\infty$$
and hence 
$$\phi(x)=\frac{2m+1}{2n+1}x+\frac{1}{2n+1}\log
\frac{(2n+1)!}{(2m+1)!}+O(e^x)+O(e^{\phi(x)}),\quad x\to -\infty,$$
which gives the first statement in \eqref{a2}.
For the second statement in (\ref{a2}) we use (\ref{formula}) and
obtain
\begin{align*}
\phi^\prime(x)&\sim \frac{(2n)!}{(2m)!}e^{(2m+1)x-(2n+1)\phi(x)}
=\frac{(2n)!}{(2m)!}e^{(2m+1)x-(2n+1)(kx+c+o(1))}\\
&=\frac{(2n)!}{(2m)!}e^{-(2n+1)c+o(1)}
=k+o(1).
\end{align*}

Let $D=\C\backslash\R_{\leq 0}$, and $p\colon D\to\C$, $p(z)=z^\mu$, the principal
branch of the power. Here $\mu$ is a complex number to be determined
so that $p$ maps $D$ onto the complement $G$ of a logarithmic spiral $\Gamma$,
with
\begin{equation}\label{6}
p(x+i0)=p(kx-i0),\quad x<0.
\end{equation}
It will be convenient to consider also the map $z\to\mu z$ obtained from $p$ by a 
logarithmic change of the variable: if $w=p(z)$ then $\log w=\mu\log z$, cf.\ Figure 1.

\begin{figure}[!htb]
\begin{center}
\captionsetup{width=.85\textwidth}
\begin{tikzpicture}[scale=0.65,>=latex](-10,-8)(10,8)
\definecolor{lightgrey}{gray}{0.93}
\filldraw [lightgrey] (-10,4) rectangle (-2,6);
\filldraw [lightgrey] (2,2.9754) -- (10,5.0246) -- (10,7.0246) -- (2,4.9754);
\draw[->](-0.7,5)->(0.7,5); 
\node at (0,5)[above]{$z\mapsto \mu z$};
\draw[->](-0.7,-5)->(0.7,-5); 
\node at (0,-5)[above]{$z\mapsto z^\mu$};
\draw[-](-10,5)->(-2,5);
\draw[-](-10,-5)->(-2,-5);
\draw[-](2,5)->(10,5);
\draw[-](2,-5)->(10,-5);
\draw[-](6,2)->(6,8);
\draw[-](-6,2)->(-6,8);
\draw[-](-6,-2)->(-6,-8);
\draw[-](6,-2)->(6,-8);
\draw[-](-10,6)->(-2,6);
\draw[-](-10,4)->(-2,4);
\draw[-](2,2.9754)->(10,5.0246);
\draw[-](2,4.9754)->(10,7.0246);
\draw[thick,-](-10,-5)->(-6,-5);
\draw[->](-6,0.7)->(-6,-0.7);
\node at (-6,0)[right]{$\exp$};
\draw[->](6,0.7)->(6,-0.7);
\node at (6,0)[right]{$\exp$};
\filldraw [black] (-6,6) circle (0.1) node[above left]{$i\pi$};
\filldraw [black] (-6,4) circle (0.1) node[below left]{$-i\pi$};
\filldraw [black] (6,6) circle (0.1) node[below right]{$i\pi$};
\filldraw [black] (6,4) circle (0.1) node[below right]{$-i\pi$};
\filldraw [black] (-3.488,6) circle (0.1);
\node at (-3,6) [above]{$a_+\!=\!x\!+\!i\pi$};
\filldraw [black] (-4,4) circle (0.1);
\node at (-2.5,4) [below]{$a_-\!=\!x\!+\log k\!-\!i\pi$};
\filldraw [black] (8.11,4.543) circle (0.1) node[below  right]{$\mu a_-$};
\filldraw [black] (8.11,6.543) circle (0.1) node[below right]{$\mu a_+$};
\draw[-,dashed](8.11,4.54)->(8.11,6.54);
\draw[-,dashed](-3.488,6)->(-4,4);
\draw [black, dotted,  thick,  domain=-1000:135, samples=200]
plot ({6-1.55*exp(0.0068*\x)*cos(\x+96.36)},{-5-1.55*exp(0.0068*\x)*sin(\x+96.36)} );
\draw [black, thick, domain=-1000:135, samples=200]
plot ({6-1.55*exp(0.0068*\x)*cos(\x+276.36)},{-5-1.55*exp(0.0068*\x)*sin(\x+276.36)} );
\draw [black, thick, dashed,  domain=-180:180, samples=100]
plot ({-6+1.6*exp(0.00447*\x)*cos(\x)},{-5+1.6*exp(0.00447*\x)*sin(\x)} );
\filldraw [black] (-9.47,-5)  arc(0:180:0.1); 
 \node at (-9.3,-5) [above left]{$-e^x$};
\filldraw [black] (-6.815,-5) arc(180:360:0.1) node[below left]{$-ke^x$};
\filldraw [black] (-4.4,-5) circle (0.1) node[below right]{$s\!=\!r^{1/\re\mu}$};
\node at (-2.79,-5.85) [below]{$=\!\sqrt{k}e^x$};
\filldraw [black] (6.2,-6.51) circle (0.1) node[below right]{$s^\mu = e^{\mu a_{\pm}}$};
\filldraw [black] (7.55,-5) circle (0.1) node[above right]{$r$};
\draw [black,thick,dashed] (6,-5) circle (1.55);
\node at (9.2,-2) [below]{$\Gamma$};
\node at (2.7,-7) [below]{$\Gamma'$};
\end{tikzpicture}
\caption{Sketch of the map $p$ and the logarithmic change of variable,
for  $k=\frac15$ and
$\mu\approx 0.9384+ 0.2403i$.
(The actual spirals $\Gamma$ and $\Gamma'$ wind much slower than drawn.)}
\label{argprinz}
\end{center}
\end{figure}
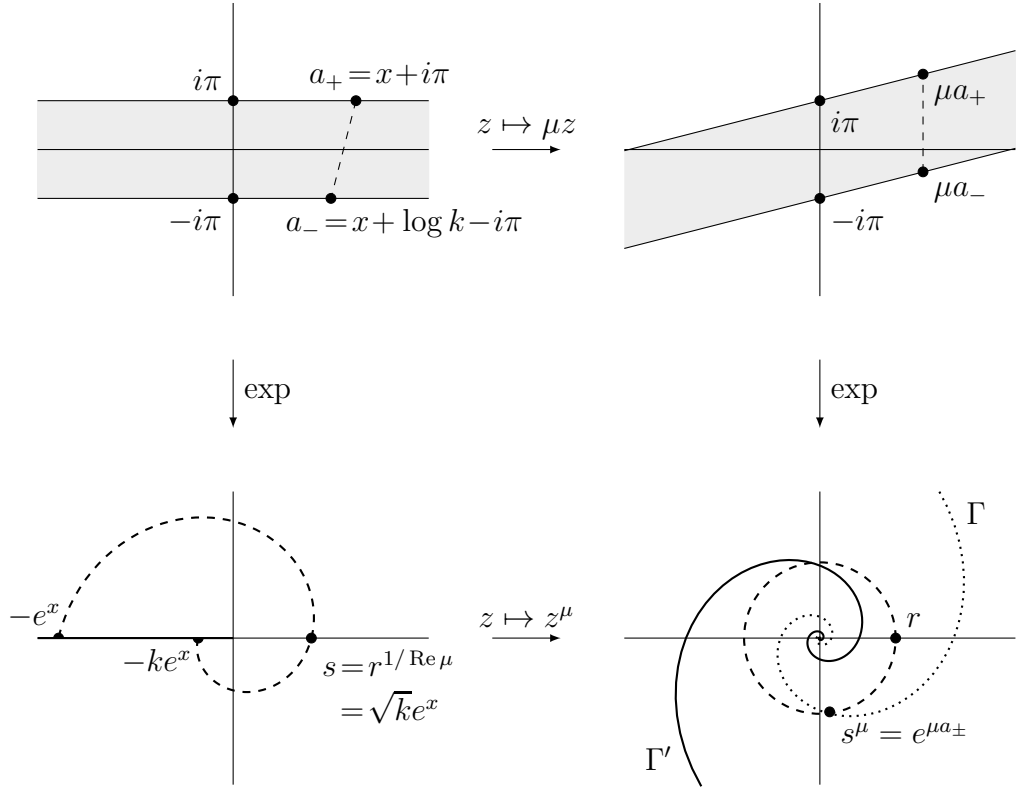

This shows (taking $x=0$ in Figure~\ref{argprinz}) that  with $a_-=\log k-i\pi$ and $a_+=i\pi$ we have
$\re(\mu a_-)=\re(\mu a_+)$; that is,
$\re(\mu (\log k -i\pi))=\re(\mu i\pi)$. Moreover, $\im(i\pi/\mu)=\pi$.
A simple computation now yields that
$$\mu=\frac{2\pi}{4\pi^2+\log^2k}(2\pi-i\log k).$$
The inverse map $h=p^{-1}$  is a conformal homeomorphism $h\colon G\to D$.
Let $\Gamma^\prime=p(\R_{\geq 0})$.
The two logarithmic spirals $\Gamma$ and $\Gamma^\prime$ divide the plane into two parts,
$G^+$ and $G^-$ which are images under $p$  of the upper and lower half-planes, respectively.

The function $V$ defined by
$$V(z)=\left\{\begin{array}{ll} 
(g_m\circ h)(z),& z\in G^+,\\
(g_n\circ h)(z),& z\in G^-,
\end{array}\right.$$
is analytic in $G^+\cup G^-$ and has a jump discontinuity
on $\Gamma$ and $\Gamma^\prime$.
In view of (\ref{a1}), (\ref{a2}) and (\ref{6}),
 this discontinuity can be removed by a small change in the independent variable.
In order to do so, we consider the strip $\Pi=\{ z\colon|\Ima z|<1\}$ and define a quasiconformal homeomorphism
$\tau\colon\C\to\C$, commuting with the complex conjugation, which is the identity outside of $\Pi$ and 
satisfies
\begin{equation}\label{Phi}
\tau(x)=\phi(x),\quad x>0,\quad\mbox{and}\quad \tau(kx)=\phi(x),\quad x<0.
\end{equation}
Our homeomorphism can be given 
by an explicit formula: for $y=\im z\in(-1,1)$ we put
$$\tau(x+iy)=\left\{ \begin{array}{ll}\phi(x)+|y|(x-\phi(x))+iy,& x\geq 0\\
\phi(x/k)+|y|(x-\phi(x/k))+iy,& x<0.\end{array}\right.
$$
The Jacobian matrix $D_\tau$ of $\tau$ is given for $x> 0$ and $0<|y|<1$ by
$$D_\tau(x+iy)=\begin{pmatrix}\phi'(x)+|y|(1-\phi'(x))&\pm(x-\phi(x))\\
0&1\end{pmatrix},$$
and we see using (\ref{a1}) that 
$$D_\tau(x+iy)\to \begin{pmatrix} 1 & 0 \\ 0 & 1 \end{pmatrix},\quad 0<|y|<1 ,\;x\to\infty,$$
Similarly, using (\ref{a2}) we find that
$$D_\tau(x+iy)\to\begin{pmatrix}1& \mp c\\ 0&1\end{pmatrix},\quad 0<|y|<1 ,\;x\to-\infty.$$
We conclude that $\tau$ is quasiconformal in the plane.

Now we modify $V$ to obtain a continuous function and define $U\colon\C\to\C$,
\begin{equation}\label{U}
U(z)=\left\{\begin{array}{ll}
(g_m\circ h)(z),& z\in G^+\cup \Gamma\cup\Gamma'\cup\{0\},\\
(g_n\circ \tau\circ h)(z),& z\in G^-.
\end{array}\right.
\end{equation}
It follows from (\ref{6}) and (\ref{Phi}) that $U$ is continuous and quasiregular in the plane.
The existence theorem for solutions of the Beltrami equation~\cite[\S V.1]{LV}
yields that there exists a a quasiconformal homeomorphism $\psi\colon\C\to\C$ with the same
Beltrami coefficient as~$U$.
The function $F=U\circ\psi^{-1}$ is then entire.

We note that $U$ is regular in $\C\backslash X$,
where $X=p(\Pi^-)$, and $\Pi^-$ is the lower half of~$\Pi$.
Let $\Delta=\{z\colon |z|>1\}$.
It is easy to see that $X\cap \Delta$ has finite logarithmic area; that is,
$$\int_{X\cap \Delta}\frac{dx\,dy}{x^2+y^2}
=\int_{\Pi^-\cap \Delta}\frac{|p'(z)|^2}{|p(z)|^2}dx\,dy
=|\mu|^2\int_{\Pi^-\cap \Delta}\frac{dx\,dy}{x^2+y^2}<\infty.$$
Thus the Beltrami coefficient of $U$  (and hence of $\psi$)
satisfies the hypotheses of the Teichm\"uller--Wittich--Belinskii theorem \cite[\S V.6]{LV}. This theorem
shows that $\psi$ is conformal at $\infty$ and may thus be normalized to satisfy
\begin{equation}\label{sim}
\psi(z)\sim z,\quad z\to\infty.
\end{equation}

Now we want to differentiate the asymptotic relation (\ref{sim}). We write  $\psi(z)=z+\psi_0(z)$
so that $\psi'(z)=1+\psi_0^\prime(z)$.
Then $|\psi_0(z)|\leq \alpha(z)$ for some function $\alpha$ satisfying $\alpha(z)=o(z)$ as $z\to\infty$.
We may assume that $\alpha(z)\to\infty$
as $z\to\infty$.
We use
the Cauchy formula
$$\psi_0^\prime(z)=\frac{1}{2\pi i}\int_{C_z}\frac{\psi_0(\zeta)}{(\zeta-z)^2}d\zeta$$
with a circle $C_z$ centered at $z$.
Choosing the radius $\beta(z)$ of this circle 
to satisfy
$$\alpha(z)=o(\beta(z)),\quad \beta(z)=o(z),\quad z\to\infty$$
and putting $Y=\{ z\colon \dist(z,X)\leq\beta(z)\}$
we obtain
\begin{equation}\label{epsilon}
\psi_0^\prime(z)\to 0,\quad z\to\infty,\;z\in\C\backslash Y.
\end{equation}
We also have
$$\meas \{\theta\in[0,2\pi]\colon re^{i\theta}\in Y\}\to 0,\quad r\to\infty.$$
Let $Y'=\psi(Y)$. Using \eqref{sim} we see that also
\begin{equation}\label{property}
\meas \{\theta\in[0,2\pi]\colon re^{i\theta}\in Y^\prime\}\to 0,\quad r\to\infty.
\end{equation}

We put
$E=F/F^\prime$.
As $F'(z)\neq 0$ for all  $z\in\C$ by construction, $E$ is entire. As all zeros of $F$ are simple, all
residues of $F'/F$ are equal to~$1$, so $E'(z)=1$ at every zero $z$ of~$E$, which implies the Bank--Laine
property. 

First we prove that $E$ is of finite order. In order to do this, we use the standard terminology
of Nevanlinna theory; see \cite{GO} or \cite{Laine}.
The counting function of the sequence
of zeros of $g_m$ and $g_n$ is of order~$1$, so the counting function of the zeros of $U$ in (\ref{U})
is of finite order. Then (\ref{sim}) shows that the counting function of zeros of~$F$, and hence the 
counting function of the zeros of~$E$,
is also of finite order; that is, $\log N(r,1/E)=O(\log r)$.
Similarly, $\log\log m(r,F)=O(\log r)$, so by the Lemma on the logarithmic derivative
\cite[Chapter~3, Theorem~1.3]{GO} we have
$\log m(r,1/E)=\log m(r,F'/F)=O(\log r)$. Thus $\log T(r,E)=O(\log r)$ so that $E$ is of finite order.

Now we estimate more precisely the growth of the Nevanlinna proximity function
$m(r,1/E)=m(r,F'/F)$. The ``small arcs lemma'' of Edrei and Fuchs \cite[Chapter~1, Theorem~7.3]{GO} permits us to
discard the exceptional set $Y'=\psi(Y)$. Outside of this set we have $\psi'(z)\to 1$ in view
of (\ref{epsilon}), therefore
\begin{equation}\label{y}
\int_{\{\theta\in[0,2\pi]\colon re^{i\theta}\in\C\backslash Y^\prime\}}\left|\log|\psi'(re^{i\theta})|\right| d\theta=o(1),
\quad r\to\infty.
\end{equation}
Furthermore, as $h(z)=z^{1/\mu}$, we have
\begin{equation}\label{x}
\int_0^{2\pi}\left|\log|h'(re^{i\theta})|\right|d\theta=O(\log r),\quad r\to\infty.
\end{equation}
Now we have in $\psi^{-1}(D^+\backslash Y)$
\begin{equation}\label{z}
\frac{F'}{F}=\left(\frac{g_m^\prime}{g_m}\circ h\circ\psi^{-1}\right)(h'\circ\psi^{-1})(\psi^{-1})'.
\end{equation}
According to (\ref{y}) and (\ref{x}), the contribution of $h'$ and $(\psi^{-1})'$
to $m(r,F'/F)$ is $O(\log r)$. 
Using the explicit form of $g_m^\prime/g_m$ we obtain, outside small neighborhoods of the zeros of $g_m$
whose contribution can be neglected again by the small arcs lemma
of Edrei and Fuchs,
\begin{equation}\label{x1}
\log^+\left|\frac{g_m^\prime(z)}{g_m(z)}\right|\sim \re^+z, \quad z\to\infty.
\end{equation}
Now the image of the circle $\{z\colon |z|=r\}$ under $h(z)=z^{1/\mu}$ is the part of the logarithmic spiral 
which connects two points on the negative real axis
and intersects the positive real axis at $r^{1/\re\mu}$; cf. Figure~1.
By \eqref{sim}, the image of this circle under $h\circ\psi^{-1}$ is an arc close to this part of the logarithmic spiral.
It now follows from \eqref{y}, \eqref{x}, \eqref{z} and \eqref{x1} that 
the part of $m(r,F'/F)$ which comes from  $\psi^{-1}(G^+\backslash Y)$
has order 
$$\rho=\frac{1}{\re\mu}=1+\frac{\log^2 k}{4\pi^2}.$$
The other part which comes from  $\psi^{-1}(G^-\backslash Y)$ is similar,
and the contribution of $Y'$ is negligible in view of (\ref{property}).
So $m(r,1/E)=m(r,F'/F)$ has order $\rho$. 

Now \eqref{B} says that
$$4A=-2\frac{E^{\prime\prime}}{E}+\left(\frac{E^\prime}{E}\right)^2-\frac{1}{E^2}.$$
It follows from the lemma on the logarithmic derivative that
$$m(r,A)=2m\left(r,\frac 1 E\right)+O(\log r).$$
Thus $A$ also has order $\rho$.

\section{Remarks}
\begin{remark}\label{rem1}
To prove that $\rho(A)=\rho$ it was sufficient to determine the growth of $m(r,1/E)$.
To show that $\rho(E)=\rho$ we also have to estimate the counting function of the zeros of~$E$.
In order to do so we note that $N(r,1/g_m)=O(r)$ and $N(r,1/g_n)=O(r)$.
Hence $N(r,1/U)=O(r^\rho)$ and thus \eqref{sim} implies that
$$N\left(r,\frac{1}{E}\right)=N(r,F)=O(r^\rho).$$
Altogether we see that $\rho(E)=\rho=\rho(A)$, as stated in the introduction.

We note that $\rho(A)<1$ implies that $\rho(E)>1$,
as follows from any of the following inequalities \cite[Theorem 12.3.1]{L}:
$$\rho(A)+\rho(E)\geq 2,\quad
\frac{1}{\rho(A)}+\frac{1}{\rho(E)}\leq 2
\quad\text{and}\quad
\rho(A)\rho(E)\geq 1.$$
Moreover, it can be deduced
from \eqref{B} that if $\rho(A)<1$,
then $\lambda(E)=\rho(E)$; see \cite[p.~442]{L}.

As our method yields examples with $\rho(E)=\rho(A)$,
it does not seem suitable to give examples with
$\rho(A)<1$.
The question whether $\rho(A)\in(\frac12,1)$ implies that
$\max\{\lambda(w_1),\lambda(w_2)\}=\infty$ for linearly independent solutions $w_1$ and $w_2$ of \eqref{DE} remains open.
\end{remark}

\begin{remark}\label{rem2}
We started our construction with two periodic locally univalent functions $g_m$
and $g_n$ and
obtained a set of orders $\rho$ which is dense in $[1,+\infty)$. By using
almost periodic
building blocks instead of $g_m$ and $g_n$, one can probably
achieve any prescribed order greater than $1$; cf. \cite[Chapter~7, Section~6]{GO}.
In this case 
$g_m$ and $g_n$ will not be explicitly known, but their asymptotic behavior
can be obtained.
\end{remark}

\begin{remark}\label{rem3}
The Bank--Laine functions we have constructed actually satisfy $E(z)=1$ whenever $E'(z)=0$.
Equivalently, one  of the two solutions of \eqref{DE} whose product is
$E$ has no zeros while the 
other one  has a finite exponent of convergence.
\end{remark}

{\em Mathematisches Seminar

Christian-Albrechts-Universit\"at zu Kiel

Ludewig--Meyn--Str. 4

24098 Kiel

Germany}
\vspace{.1in}

{\em Purdue University

West Lafayette, IN 47907

USA}
\end{document}